# ON THE CAPILLARITY EQUATION IN TWO DIMENSIONS

Rajat Bhatnagar and Robert Finn

**Summary:** We study the capillarity equation from the global point of view of behavior of its solutions without explicit regard to boundary conditions. We show its solutions to be constrained in ways, that have till now not been characterized in literature known to us.

The capillarity equation

$$\operatorname{div} Tu = \kappa u, \quad Tu \equiv \nabla u / \sqrt{1+|\nabla u|^2} \tag{1}$$

describes the height $u(x,y)$ of a surface whose mean curvature $H$ is proportional to its height, with factor $\kappa/2$. Physically, $u(x,y)$ represents an (ideally thin) interface separating two fluids, with a density difference $\rho$ (heavier fluid below when $\kappa > 0$), in a uniform gravity field $g$ directed downward. Here $\kappa = \rho g / \sigma$ is the *capillarity constant*, with $\sigma$ the surface tension of the fluid interface.

In a typical physical problem, the surface is to be defined over a horizontal section of a cylindrical region $\Omega$ with vertical bounding walls, which it is required to meet in a prescribed angle $\gamma$. Historically, the most familiar example is the circular glass "capillary tube" of small radius, which when dipped into water exhibits spectacular rise of the liquid in its interior. In the citations [1 - 5], $\Omega$ is chosen to be a section between two vertical and parallel plates of infinite extent, partially immersed into an unlimited fluid mass, that is horizontal and at rest at infinite distance from the plates. This configuration has the virtue that all solutions of (1) between the plates, and meeting them in two constant (perhaps distinct) contact angles, can be expressed in terms of the single coordinate orthogonal to both plates (see [1], where this configuration was first studied intensively, and [6] from which that assertion can be proved). Denoting the indicated coordinate by $x$, the equation (1) simplifies to the two dimensional case:

$$\left( \frac{u_x}{\sqrt{1+u_x^2}} \right)_x \equiv (\sin \psi)_x = \kappa u \tag{2a}$$

where $\psi$ is inclination angle of the curve $u(x)$, so that

$$u_x = \tan \psi \,. \tag{2b}$$

Geometrically, (2a) asserts that the planar curvature of the solution curve is proportional to its height $u$; this is the essential substance underlying much of the material in the indicated citations above. This comment has also striking consequences for global study of the equation, without specific regard to boundary conditions. We focus on that perspective in what follows.

## 1. General considerations; attracting solutions.

We observe initially that at points where $u \neq 0$, (2a,b) is equivalent to the first order system

$$\frac{dx}{d\psi} = \frac{\cos \psi}{\kappa u}, \quad \frac{du}{d\psi} = \frac{\sin \psi}{\kappa u}. \tag{3a,b}$$

and admits a first integral of the form

$$\frac{\kappa}{2}u^2 + \cos\psi = c \equiv const. \tag{4}$$

see, e.g., Gyamant [14], or alternatively integrate (3b) which separates.

The integral (4) can be extended to a complete integral via the relations

$$c^* = \frac{\kappa}{2}(u^*)^2 + \cos\psi^*$$

$$x = x^* + \frac{1}{\sqrt{2\kappa}} \int_{\psi^*}^{\psi} \frac{\cos\tau\, d\tau}{\sqrt{c^* - \cos\tau}}, \tag{5}$$

$$u = \sqrt{\frac{2}{\kappa}}\sqrt{c^* - \cos\psi}$$

determining $x$ and $u$ in terms of prescribed data $x^*$ and $u^*$ at a chosen $\psi^*$. We see immediately from (5) the translation invariance of solutions, in the variable $x$, and also the reflection invariance in the $x$ and $u$ coordinates.

In what follows, we make the useful observation that in terms of the non-dimensional coordinates $\xi = \sqrt{\kappa}x, U = \sqrt{\kappa}u$, our relations above simplify to universal relations with no parameters:

$$(\sin\psi)_\xi = U \qquad U_\xi = \tan\psi . \tag{2*a,b}$$

$$\frac{d\xi}{d\psi} = \frac{\cos\psi}{U}, \quad \frac{dU}{d\psi} = \frac{\sin\psi}{U} \tag{3*a,b}$$

$$\frac{1}{2}U^2 + \cos\psi = c \equiv const. \tag{4*}$$

thus facilitating physical discussion and interpretation; similarly (5) changes to:

$$C^* = \frac{1}{2}(U^*)^2 + \cos\psi^*$$

$$\xi = \xi^* + \frac{1}{\sqrt{2}} \int_{\psi^*}^{\psi} \frac{\cos\tau\, d\tau}{\sqrt{C^* - \cos\tau}} \tag{5*}$$

$$U = \sqrt{2}\sqrt{C^* - \cos\psi}$$

Note that in accord with differing topics of the two papers, the present choice of non-dimensional coordinates differs from the choice appearing in [12].

For the purpose of constructing solutions of (3*a,b) as graphs $U(\xi)$ over the $\xi$ – axis, we limit ourselves to the range $-\pi/2 < \psi < \pi/2$. Using (5) rewritten as (5*) in the universal coordinates, with $U^* = 0$, $\xi^* = 0$, $\psi^* = \psi_0$, we find that the apparent singularity at $U = 0$ is to some extent deceptive, and that solutions extending continuously across the $\xi$ – axis with prescribed angle of crossing are (uniquely) obtainable explicitly from (5*). However, for the solutions thus constructed, the orientation of the solution curve reverses as the axis is crossed, see Figure 1; thus, for curves crossing the axis, the integration in (5*) must be split into two separate terms.

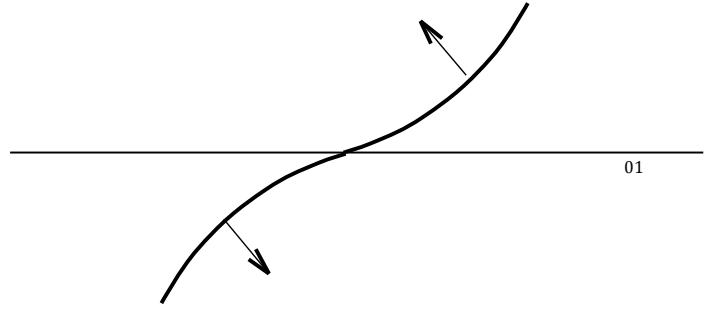

**Figure 1.** Reversal of curvature orientation on crossing *x* – axis; the sense of increasing $\psi$ reverses.

A more unified (and intrinsic) interpretation appears by introducing arc length *s* as parameter for describing the curve; we obtain the system

$$\frac{d\xi}{ds}=\sqrt{1-v^2},\ \frac{dU}{ds}=v,\ \frac{dv}{ds}=U\sqrt{1-v^2} \qquad (3^{**})$$

with $v = \sin\psi$. Here the description is intrinsic, with no need to take account of orientation reversal, and one sees directly the global analyticity of the solutions in *s*, on any interval on which |v| < 1, see, e.g., [13], pp. 32-37. This latter condition is automatically imposed by the limitation to graphs, for all configurations arising in this paper. We remark however that a branching does occur, whenever the value *v* = ±1 appears at any point. It is then effected by a choice of sign for the square root, in the continuation through that point.

We use here the same symbol *U* in a generic sense for the various functions $U(\xi)$, $U(\psi)$, $U(s)$, $U(x)$.

The latter two equations of $(3^{**})$ lead to the equivalent second order equation for the height *U*:

$$U_{ss}=U\sqrt{1-U_s^2} \qquad (6)$$

in terms of arc length *s* as independent variable. This relation admits the first integral $U^2 + \sqrt{1-v^2} = C$, on which the theory could alternatively have been based. However in what follows, we have found the traditional parameter $\psi$ to be more convenient for the exposition.

Much information can be gleaned directly from $(4^*)$. In the indicated range for $\psi$ there can be no solution curve for which $c \leq 0$. If $c > 1$, then solutions $U(x)$ cannot change sign continuously, and they are bounded below in magnitude. Setting $c - 1 = \delta$, we obtain $U^2 > 2\delta > 0$, and thus by $(2^*a)$ the planar curvature *k* satisfies $k^2 > 2\delta > 0$. We are led to the general statement:

**Theorem 1.** *If $c > 1$ in $(4^*)$, then the length l of any $\xi$ interval over which the solution curve can extend satisfies $l < \sqrt{2/\delta}$. Any solution curve positive at any point can be extended across a unique positive minimum $U_0$, and following a horizontal translation to put the minimum*

*point at $\xi = 0$, the set of all solutions simply covers a domain $\mathcal{D}^+$ as indicated in* Figure 2. *Here each curve of the family is a solution of* ($3^*$a,b) *in the parameter interval* $-\pi/2 < \psi < \pi/2$.

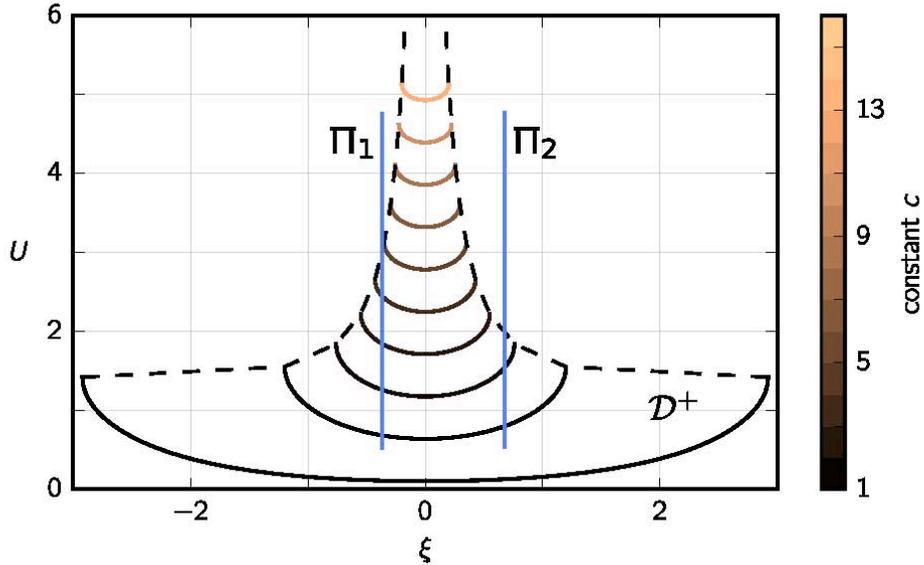

**Figure 2. The set of all positive graph solutions of $2^*$, modulo horizontal translation; each solid curve is a "curve of constant attracting force".**

The curves in Figure 2 admit an interesting physical interpretation. We consider them as vertical sections of the surface of an infinite channel of fluid orthogonal to the plane of the figure, and we imagine two vertical parallel plates, partly inserted into the fluid as indicated in the figure. We suppose that the contact angles at the triple interface on the sides of the plates that face each other are those indicated in the figure. Then there is a horizontal attracting force between the plates, relative to the surface tension $\sigma$, of magnitude $\mathcal{F} = U_0^2 = \kappa u_0^2$. *This force is independent of the particular positions of the plates along the curve*; it depends only on the height of the curve at the (single) minimizing point. Thus the indicated region is swept out by "*curves of constant attracting force*".

For a proof of the stated force formula, see, e.g., [3].

Each curve in Figure 2 can be viewed as the "canonical representative" of an equivalence class of solution curves obtained by rigid horizontal translation, and yielding the same force relation. There is additionally the freedom to reflect the entire configuration in the $\xi$–axis, yielding an equivalent example with negative heights, but the same attracting forces. In general, also reflection in a vertical axis is permissible, however in the present (symmetric) case, that would lead to no change. Subject to these remarks, *the configuration illustrated includes all conceivable attracting configurations*, in the sense that any attracting case is included as one of the illustrated curves. That is, *the set of all attracting cases is covered by the choice of values $c > 1$ in* ($4^*$).

We emphasize that the theory does not encompass changes in scale, beyond that imposed by the physical parameter $\sqrt{\kappa}$. For an air/water interface in the earth's gravity field, one finds $\kappa \approx 14\,\text{cm}^{-2}$.

**The case c = 1.** We look next for positive solutions when $c = 1$. From (4*) we find directly the explicit function

$$U = 2\sin\frac{\psi}{2} \tag{7}$$

in terms of its inclination $\psi$, from which using (3*a)

$$\frac{d\xi}{d\psi} = \left(\frac{1}{2\sin\frac{\psi}{2}} - \sin\frac{\psi}{2}\right) \tag{8}$$

so that

$$\xi = \xi_2 + \log\frac{\tan(\psi/4)}{\tan(\psi_2/4)} + 2\left(\cos\frac{\psi}{2} - \cos\frac{\psi_2}{2}\right) \tag{9}$$

together with (7) provide a solution inclined at angle $\psi_2$ when $\xi = \xi_2$, which is positive and increasing for $\xi < \xi_2$, and which vanishes at $\xi = -\infty$

We identify this curve below as the solution "**I**" of the citations [3,4]. It is defined and positive on the half-space $\xi < 0$, with $U_\xi > 0$, vanishes with its slope at $x = -\infty$, and is vertical at $\xi = \xi_3$, for a $\xi_3 > \xi_2$, see Figure 1 of [3]. Although (7,9) is achieved as a uniform limit over compact subdomains of the $\xi$– axis, of curves with horizontal symmetry, it is itself not symmetric, and it has lost the positive bound from below that occurs on each of the approximating curves. *It can be characterized up to horizontal translation and reflection in the coordinate axes, as the only solution of (2*a,b) that extends to infinity in any direction.*

Four metrically identical copies of this curve can be obtained, all of them solutions of (2*a,b), by reflection in the coordinate axes. To interpret the curve we have found in terms of Figure 2, we focus attention on the lowest of the solution curves shown in the figure, then translate that curve horizontally to situate the point at which $\psi = \psi_2$ on the line $\xi = \xi_2$. We do that for a sequence of curves corresponding to a sequence of values $U_{0j} \to 0$ (see Figure 3). *The limit curve of that procedure is uniquely determined as the curve "**I**" of the earlier papers [3,4]. It has the property that the net force on any two vertical plates that are partially submerged and cut through it orthogonal to the plane of the page is zero.*

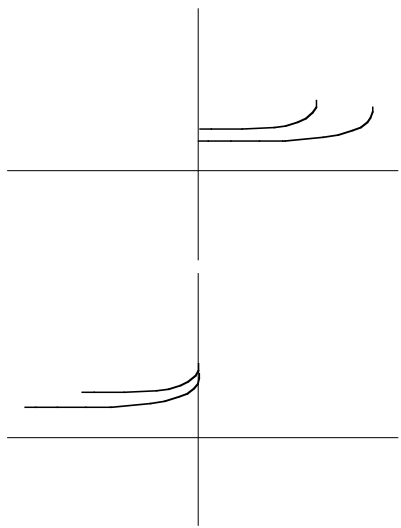

**Figure 3*a*.** Two curves of the attracting family of Figure 2, prior to and following the shifting: This is a conceptual sketch, not to scale.

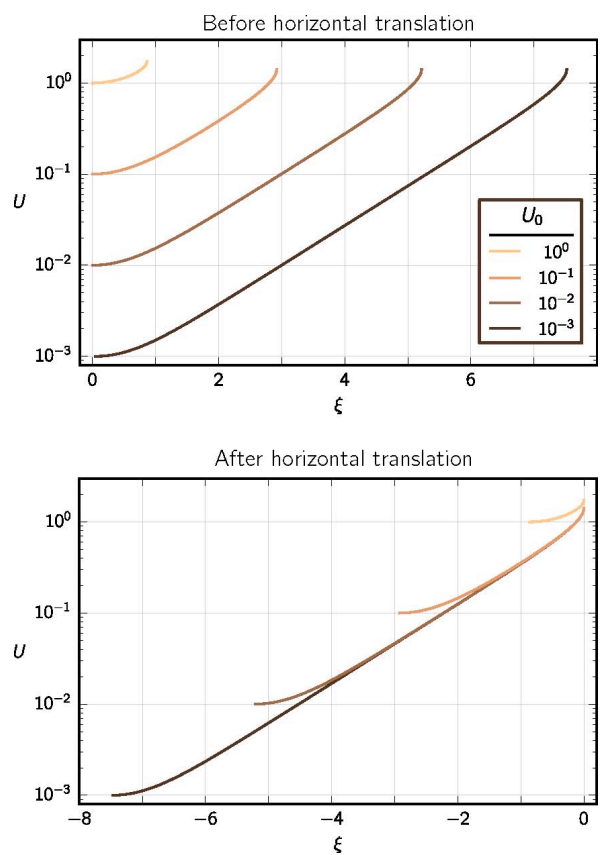

**Figure 3*b*.** Actual appearance of initial and translated solution curves.

## 2. Repelling solutions; $0 < c < 1$.

As shown in [3], the repelling solution curves joining the plates cross the $\xi$ – axis in a non-zero angle $\psi_0$, and the net normalized repelling force is

$$\mathcal{F} = -2(1 - \cos\psi_0), \tag{10}$$

the negative sign serving to distinguish it from attracting forces. Analogous to the attracting case, any two parallel vertical plates that cut a given solution curve and adopt the incident contact angles on the sides facing each other will repel each other with the force (10).

Again we may normalize the solution curves by horizontal translation, obtaining as representatives of equivalence classes, curves that cross the $\xi$ – axis at the origin of coordinates. The heights for $\xi < 0$ will then be the negatives of the corresponding heights for $\xi > 0$. The equation (4*) yields now

$$U^2 = 2(\cos\psi_0 - \cos\psi). \tag{11}$$

We consider a particular solution curve passing through the origin, and inclined at that point with angle $\psi_0$, $0 < \psi_0 < \pi/2$. Thus $\sin\psi_0 > 0$. From (2*) we see that at $\xi = 0$, $(\sin\psi)_\xi = 0$, $(\sin\psi)_{\xi\xi} = U_\xi(0) = \tan\psi_0 > 0$. Thus, the curve inflects at the origin, $(\sin\psi)_\xi$ is increasing at $\xi = 0$, and from (2*) one finds that $\sin\psi$ then continues to increase to the value $\sin\psi = 1$, beyond which further continuation as a graph is not possible. We see that the solution segment is confined to an interval $-\xi_0 < \xi < \xi_0$, at the endpoints of which the curve is vertical. From (11) we find that at the endpoints

$$U_0^2 \equiv U^2(\pm\xi_0) = 2\cos\psi_0. \tag{12}$$

To complete the characterization, we determine $\xi_0$ in terms of $\psi_0$. For a generic value $\xi$, we find by inserting (12) into (3*a)

$$\sqrt{2}\xi = \int_{\psi_0}^{\psi} \frac{\cos\tau}{\sqrt{\cos\psi_0 - \cos\tau}} d\tau \tag{13}$$

Writing $\cos\tau = -s$, (13) becomes

$$\sqrt{2}\xi = -\int_{s_0}^{s} \frac{t}{\sqrt{1-t^2}} \frac{1}{\sqrt{t-s_0}} dt \tag{14}$$

so that

$$\sqrt{2}\xi_0 = -\int_{s_0}^{0} \frac{t}{\sqrt{1-t^2}} \frac{1}{\sqrt{t-s_0}} dt. \tag{15}$$

One sees directly that as $\psi_0 \to 0$, $s_0$ decreases monotonically to $-1$, so that $\xi_0 \to \infty$. But as $\psi_0 \to \pi/2$ then $s_0 \to 0$, and we find from (15) that $\xi_0 \to 0$.

We can infer from this information the general character for large $|x|$ of the region $\mathcal{E}$ swept out by the solution curves for repelling configurations, as shown in Figure 4. To justify the figure as shown near the origin, we must calculate $\partial\xi_0/\partial U_0$ from (15), using $s_0 = -\cos\psi_0$. We cannot differentiate (15) directly under the sign, as that leads to a singular integral. But we can integrate by parts, obtaining

$$\sqrt{2}\xi_0 = 2\int_{s_0}^{0} \frac{\sqrt{t-s_0}}{(1-t^2)^{3/2}} dt. \tag{16}$$

From (12) we calculate directly, using (3*b)

$$\frac{\partial \xi_0}{\partial U_0} = \frac{\partial \xi_0}{\partial s_0}\frac{\partial s_0}{\partial \psi_0}\frac{\partial \psi_0}{\partial U_0} = -\frac{U_0}{\sqrt{2}}\int_{s_0}^{0}\frac{dt}{\left(1-t^2\right)^{3/2}\sqrt{t+\cos\psi_0}} \qquad (17)$$

which clearly vanishes in the limit as $U_0 \to 0$.

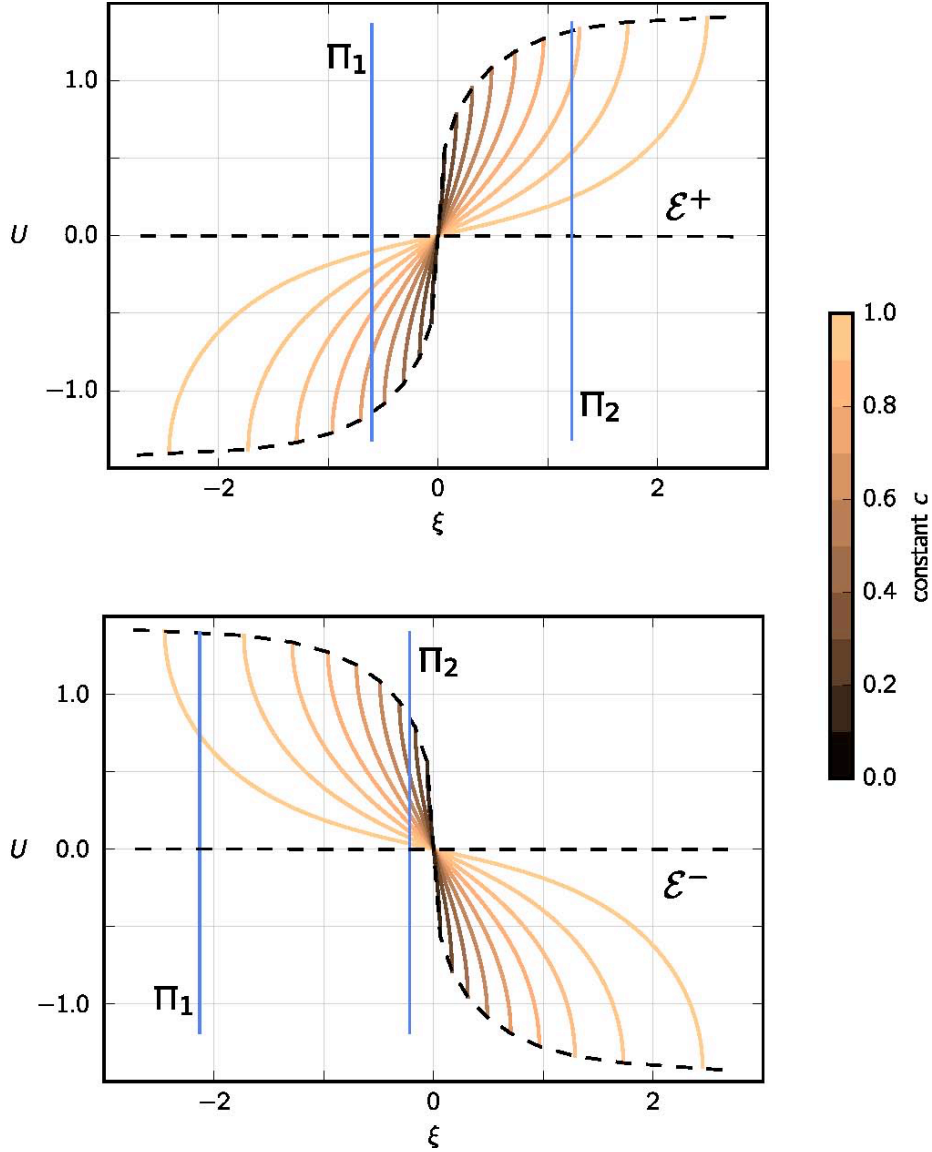

**Figure 4.** The regions $\mathcal{E}^+$ and $\mathcal{E}^-$ swept out by the solution curves in repelling case. The two families of curves are respectively distinct. The dashed curves are the loci of vertical points on solution curves. Any two vertical plates meeting a common curve $\mathcal{C}$ in angles $\gamma_1$ and $\gamma_2$ as indicated repel each other with the same force $\mathcal{F} = 2(1 - \cos\psi_0)$, where $\psi_0$ is the angle with which $\mathcal{C}$ cuts the horizontal axis.

From (13) and (15) we conclude:

**Theorem 2.** *All repelling solutions are contained between the horizontal lines $U_0^2 = \pm 1$. When normalized to pass through the origin, their vertical points determine two domains $\mathcal{E}^+$ and $\mathcal{E}^-$ as indicated in Figure 4. Any two vertical plates $\Pi_1$ and $\Pi_2$ cutting a common curve $\mathcal{C}$ in contact angles $\gamma_1$ and $\gamma_2$ as indicated in the figure yield the same normalized horizontal force $\mathcal{F}$, depending only on $\mathcal{C}$. Only forces in the range $0 \leq |\mathcal{F}| < 2$ can be obtained. Except for the*

*trivial case $U \equiv 0$, each such force can arise from exactly two curves of the indicated family, cutting the x – axis in an angle $\psi_0$ within the admissible range $-\pi/2 < \psi_0 < \pi/2$.*

If we continue the construction indicated in Figure 4, introducing further solutions and extending the dashed lines, while shifting the curves horizontally to keep new points at which $\psi = \pi/2$ centered at $\xi = 0$, we are led again by an analogous limiting procedure to the curve "**I**" of the citations [3,4], and to its reflections in horizontal and vertical axes.

**Final Remarks.** Readers may have noticed that we have omitted the important case $\kappa < 0$, which governs (for example) the base of a fluid column emerging from a "medicine dropper". In this case gravity is again downward directed, but the heavier fluid lies above the surface interface. The theory then takes a different form, leading to configurations very different from those studied here, and it links inextricably with stability criteria, which in turn are linked with the presence of rigid boundaries. In this respect we have taken the view that the topic would be out of character for the present limited discussion. There is an extensive literature for such questions, for which we cite, among many other studies, the papers [7-11], all devoted to a particular theme of special interest, arising within the theory. Those citations address only one of the many striking phenomena that can occur. There is in fact little overlap with the material offered in the present study, which seems to have been largely overlooked prior to the citations immediately preceding it.

The case $\kappa = 0$ yields formally minimal surfaces, which appear, e.g., as soap films spanning a closed wire. For a volume-constrained problem in absence of gravity, one obtains surfaces of constant mean curvature. That special case is encountered notably in problems arising in space travel, and is again outside the purview of the present study.

We remark finally that alternative derivations in a more general context, of the force formulas used here, can be found in [5] and on p. 122 of [15].

## Acknowledgments.

The second author is indebted to the Mathematische Abteilung der Universität, in Leipzig, for its hospitality and for excellent working conditions, during his visit of October 2014. He is indebted to the Simons Foundation for a very helpful grant facilitating his travel. He thanks Erich Miersemann for incisive comments on some conceptual points related to the problem considered.

Department of Biochemistry and Biophysics,
California Institute for Quantitative Biosciences,
University of California
San Francisco, CA 94143, USA
raj.bhatnagar@gmail.com

Mathematics Department
Stanford University
Stanford, CA 94025, USA
finn@math.stanford.edu